\documentclass[a4paper,11pt]{article}
\usepackage[latin1]{inputenc}
\usepackage{amssymb}
\usepackage{a4wide}
\usepackage{amsmath}
\usepackage{calc}

\newcommand{\begr}[1]{\textbf{#1}}

\newenvironment{rom}{\renewcommand{\theenumi}{(\roman{enumi})}
\renewcommand{\labelenumi}{\theenumi}
\begin{enumerate}}{\end{enumerate}}

\newenvironment{entry}{\begin{list}{}{\setlength{\labelwidth}{68pt}%
\setlength{\itemsep}{\smallskipamount}\setlength{\leftmargin}{\labelwidth+\labelsep}}}{\end{list}}

\newcommand{\CN}{\mbox{\ensuremath{\mathit{CN}}}}
\newcommand{\ON}{\mbox{\ensuremath{\mathit{ON}}}}

\newcommand{\womega}{\omega}

\begin{document}

\title{The $f$-Factor Problem for Graphs and the Hereditary Property\footnote{%
This paper was supported by the Volkswagen Stiftung}}
\author{Frank Niedermeyer, Bonn\\Saharon Shelah, Jerusalem\\Karsten Steffens, Hannover}
\date{}

\maketitle

\begin{abstract}
If $P$ is a hereditary property then we show that, for the existence of a perfect 
$f$-factor, $P$ is a sufficient condition for countable graphs and yields a sufficient
condition for graphs of size $\aleph_1$. Further we give two examples of a hereditary
property which is even necessary for the existence of a perfect $f$-factor. We also
discuss the $\aleph_2$-case.
\end{abstract}

We consider graphs $G=(V,E)$, where $V=V(G)$ is a nonempty set of vertices and 
$E=E(G)\subseteq \{$ \mbox{$e\subseteq V\colon$} \mbox{$|e| = 2$}$\}$ is the set of 
edges of $G$. If $x$ is a vertex 
of $G$ and $F \subseteq E$, then we denote by $d_F(x)$ the cardinal 
$|\{ e\in F\colon x\in e\}|$. $d_F(x)$ is called the \begr{degree of $x$ with
respect to $F$} and $d_E(x)$ the \begr{degree of $x$}. $\ON$ denotes the class
of ordinals, $\CN$ the class of cardinals. Greek letters $\alpha, \beta, \gamma, \ldots$
always denote ordinals, whereas the middle letters $\kappa, \lambda, \mu, \nu, \ldots$
are reserved for infinite cardinals.

Let $G = (V,E)$ be a graph, $f\colon V\rightarrow CN$ be a function and $F \subseteq E$.
$F$ is said to be an \begr{$f$-factor} of $G$ if $d_F(x) \leq f(x)$ for all $x\in V$.
We call an $f$-factor $F$ of $G$ \begr{perfect} if $d_F(x) = f(x)$ for all $x\in V$.
For $\kappa\in CN$ we denote $f^{-1}(\kappa):=\{x\in V:f(x)=\kappa\}$. 

Let $\cal C$ be the class of all ordered pairs $(G,f)$, such that $G = (V,E)$ is a graph,
$f\colon V\rightarrow CN$ is a function, and $f(x)\leq d_E(x)$ for all $x\in V$.

This paper discusses the problem to find a necessary and sufficient condition
for the existence of a perfect $f$-factor of a graph. In [5], Tutte published a criterion
for finite graphs, and in [4] Niedermeyer solved the problem for countable graphs and functions
$f:V\longrightarrow \omega $. We present a solution for graphs of size $\aleph_0$ and functions
$f:V\longrightarrow \omega \cup\{\aleph _0\}$, a solution for graphs of size $\aleph _1$, and
discuss the $\aleph_2$-case.

If $H \subseteq E$, then denote by $G-H$ the graph $(V,E\setminus H)$, and if
$e\in E$, then let $G-e$ be the graph $G-\{e\}$. If $x,y\in V$, denote by
$f_{x,y}\colon V\rightarrow \CN$ the function defined by $$f_{x,y}(v) := 
\begin{cases}f(v)-1 & \text{if }v\in\{x,y\}\text{ and }1\leq f(v) < \aleph_0 \\
             f(v)   & \text{else}\end{cases}.$$ 
Now let $P$ be a formula with two free variables.
$P(G,f)$ means that
$(G,f)\in\cal C$ and
$(G,f)$ has the property $P$. $P$ is said to be \begr{hereditary} if for every
$(G,f)$ with $P(G,f)$, for every vertex $x\in V(G)$ with $f(x) > 0$ there exists a vertex
$y\in V(G)$ with $f(y) > 0$, $\{x,y\} \in E(G)$, and $P(G - \{x, y\}, f_{x,y})$.

\paragraph{Remark} Let $P$ be a hereditary property, let $(G,f)\in\cal C$ such that
$P(G,f)$, and let $W\subseteq V(G)$ be finite. Then there exists a finite $f$-factor
$F$ of $G$ such that $P(G-F,f-d_F)$, $d_F(x)=f(x)$ for every
$x\in W$ with $f(x)<\aleph_0$, and $d_F(x)>0$ for every $x\in W$ with 
$f(x)\geq \aleph_0$.

\paragraph{Example~1} Let $P_1(G,f)$ be the property ``$G$ possesses a perfect
$f$-factor''. Obviously $P_1$ is a hereditary property.

\paragraph{Definition} Let $(G,f)\in\cal C$.
By recursion on $\alpha \in \ON$ we define the property that $(G, f)$ is
an \begr{$\alpha$-obstruction}. Let $G=(V,E)$.

If there is an $x\in V$ with $f(x) > 0$ such that $f(y) = 0$ for all $y\in V$ with
$\{x,y\}\in E$, then $(G, f)$ is a \begr{$0$-obstruction}.

If there is a vertex $x\in V$ such that $f(x) > 0$ and
\begin{rom}
\item for every $y\in V$ with $\{x,y\}\in E$ and $f(y)>0$ there is an 
ordinal  $\beta_y$
such that $(\mbox{$G-\{x,y\}$}, f_{x,y})$ is a $\beta_y$-obstruction and
\item $\alpha = \sup \{\beta_y + 1\colon \{x,y\}\in E, f(y) > 0 \}$,
\end{rom}
then $(G, f)$ is an \begr{$\alpha$-obstruction}.

\paragraph{Example~2} Let $P_2(G, f)$ be the property ``$(G, f)$ is not an $\alpha$-obstruction
for every $\alpha\in\ON$''. Then we can prove the
following

\paragraph{Lemma~1}
\begin{rom} 
\item $P_2$ is a hereditary property.
\item If $P$ is a hereditary property, then $P_2(G,f)$ is necessary for $P(G,f)$.
Therefore $P_2$ is a necessary condition for the existence of a perfect $f$-factor.
\end{rom}

\paragraph{Proof}
\begin{rom}
\item Assume $P_2(G, f)$, that means that for all $\alpha \in \ON$, $(G,f)$ is
not an $\alpha$-obstruction. Let $G=(V,E)$ and $x\in V$ with $f(x)>0$. To get a
contradiction let us assume that, for each $y\in V$ with $\{x,y\} \in E$ and $f(y)>0$,
there is an ordinal $\beta_y$ such that $(G-\{x,y\}, f_{x,y})$ is a 
$\beta_y$-obstruction. If $\alpha = \sup \{\beta_y + 1\colon \{x,y\}\in E, f(y)>0\}$, 
then $(G,f)$ is an $\alpha$-obstruction which contradicts our assumption.

\item By induction on $\alpha\in\ON$ we prove for any $(G,f)\in\cal C$ with $P(G,f)$
that $(G,f)$ is not an $\alpha$-obstruction.

Since $P$ is heriditary, $(G,f)$ is obviously not a $0$-obstruction.

Now let $\alpha > 0$. Assume that $(G,f)$ is an $\alpha$-obstruction. Let $G=(V,E)$.
By definition, there is a vertex $x\in V$ with $f(x) > 0$ such that for each 
$y\in V$ with $f(y)>0$ and $\{x,y\}\in E$ there is an ordinal $\beta_y<\alpha$ such that
$(G-\{x,y\}, f_{x,y})$ is a $\beta_y$-obstruction. On the other hand, since $P(G,f)$, $P$ is
hereditary, and $f(x) > 0$, there is an edge $\{x,y\}\in E$ such that $P(G-\{x,y\},f_{x,y})$.
By inductive hypothesis $(G-\{x,y\}, f_{x,y})$ is {\it not} a $\beta_y$-obstruction. This
contradiction proves (ii).
\end{rom}

\bigskip

For a hereditary property $P$, it must not be true that $P_2(G,f)$ is sufficient for
$P(G,f)$. This is demonstrated by the following example.

\paragraph{Example 3} Let $P_3(G,f)$ be the property ``$G$ possesses a perfect $f$-factor
without cycles''.

\bigskip

$P_3$ also shows that not every hereditary property is a necessary condition for the
existence of a perfect $f$-factor.

\paragraph{Definition} Let $(G,f)\in\cal C$.
For $0<k\leq \womega$ we call a sequence $T=(v_i)_{0\leq i<k}$ of vertices of $G$
a \begr{trail} if $\{v_{i-1}, v_i\} \in E(G)$ for $0<i<k$ and $\{v_{i-1}, v_i\} \neq
\{v_{j-1}, v_j\}$ for $i\neq j$. For any $f$-factor $F$, a trail $T=(v_i)_{0\leq i<k}$
is called \begr{$F$-augmenting} if
\begin{rom}
\item $k>1$
\item $\{v_{i-1}, v_i\} \in F$ iff $i>0$ is even
\item $d_F(v_0) < f(v_0)$
\item $k = \womega$ \\ 
      \hbox{}\hspace{1ex}or \\
      $k < \womega$ is even, $v_0 \neq v_{k-1}$ and $d_F(v_{k-1}) < f(v_{k-1})$\\
      \hbox{}\hspace{1ex}or \\
      $k < \womega$ is even, $v_0 = v_{k-1}$ and $d_F(v_{k-1})+1 < f(v_{k-1})$
\end{rom}

\paragraph{Example~4} Let $P_4(G,f)$ be the property
``for every $f$-factor $F$ of $G$ and every vertex $x\in V(G)$ with
$d_F(x) < f(x)$ there exists an $F$-augmenting trail starting at $x$''.
Further let $P_4'(G,f)$ be the property ``$P_4(G,f)$ and
$ran(f)\subseteq\womega$''.

\paragraph{Lemma~2} If $(G,f)\in\cal C$
and $G$ possesses a perfect $f$-factor, then $P_4(G,f)$.

\paragraph{Proof}
For the convenience of the reader, we present the easy proof. Let $G=(V,E)$, let $F$ be an
$f$-factor of $G$ and $H$ be a perfect $f$-factor of $G$. For all $x\in V$ with $d_F(x) < f(x)$,
we construct by induction an $F$-augmenting trail starting at $x$. Let $v_0 = x$.
Since $d_F(v_0) < f(v_0) = d_H(v_0)$ there is an edge $\{v_0,y\}\in H\setminus F$.
Let $v_1 = y$. Let the trail $T=(v_j)_{0\leq j\leq i}$
be defined such that 
\begin{itemize}
\item[(1)] $\{v_{j-1}, v_j\} \in F\setminus H$ iff $j>0$ is even.
\item[(2)] $\{v_{j-1}, v_j\} \in H\setminus F$ iff $j$ is odd.
\end{itemize}
If $i$ is odd, $v_i\neq v_0$, and $d_F(v_i) < f(v_i)$, let $k=i+1$.

If $i$ is odd, $v_i = v_0$, and $d_F(v_i)+1 < f(v_i)$, let again $k=i+1$.

If $i$ is odd and
$v_i \neq v_0$, $d_F(v_i) = f(v_i)$ or $v_i=v_0$, $d_F(v_i)+1\geq f(v_i)$, then there
is an edge $\{v_i,y\}\in F\setminus H$ which
is not an edge of $T$. Let $v_{i+1} = y$.

Finally, if $i$ is even, there is an edge $\{v_i, y\}\in H\setminus F$
which is not an edge of the trail $T$. Let $v_{i+1} = y$.

\bigskip

Much more difficult is the proof of Lemma~3 which is Corollary~4 of \cite{4}.

\paragraph{Lemma~3} $P_4'$ is a hereditary property.
\bigskip

It is not true that every hereditary property $P$ is a sufficient condition for the
existence of a perfect $f$-factor of a given graph. This demonstrates the property
$P_4$, applied to the complete bipartite graph $K_{\aleph_0, \aleph_1}$ and the
function $f\equiv 1$. But we have the following

\paragraph{Theorem~1} Let $(G,f)\in\cal C$ and $|V(G)| = \aleph_0$. If $P$ is a hereditary
property and $P(G,f)$ then $G$ possesses a perfect $f$-factor.

\paragraph{Proof}
Let $v_0, v_1, v_2, \ldots$ be an enumeration of the vertices of $G$ such that,
for every $x\in V$ with $f(x) = \aleph_0$, the set $\{i<\womega\colon x=v_i\}$ is 
infinite.  Since $P(G,f)$ and $P$ is hereditary, one can define recursively finite 
$f$-factors $F_0 \subseteq F_1 \subseteq F_2 \subseteq \cdots $ such that 
$(G-F_k, f-d_{F_k})$ fulfills property $P$ and the following is true: If 
$f(v_0) = \aleph_0$, then $F_0$=$\{\{x,v_0\}\}$, if $f(v_k)=\aleph _0$, $k>0$, then  
$F_k\setminus F_{k-1} = \big\{\{x, v_k\}\big\}$ for some $x\in V$, and if
$f(v_k) < \aleph_0$, then $d_{F_k}(v_k) = f(v_k)$. By construction, 
$F:= \bigcup \{ F_k \colon k < \womega \}$ is a perfect $f$-factor.

\paragraph{Corollary~1}
Let $(G,f)\in\cal C$ and $|V(G)| = \aleph_0$. 
\begin{itemize}
\item[(1)] $G$ has a perfect $f$-factor iff $P_2(G, f)$.
\item[(2)] If $ran(f)\subseteq\womega$, then $G$ has a perfect $f$-factor iff $P_4(G,f)$.
\end{itemize}
Tutte's condition $(\cite{3},\cite{5})$ for the existence of a perfect $1$-factor for finite
graphs is necessary but not sufficient for countable graphs. Thus Theorem~1 shows that
not every necessary condition for the existence of a perfect $f$-factor is a hereditary
property. The property ``$G$ has a perfect $f$-factor
with cycles'' tells us that a sufficient condition for the existence of a perfect
$f$-factor for $G$ is not necessarily hereditary.

\paragraph{Definition} Let $(G,f)\in\cal C$, $G=(V,E)$, and $|V|=\kappa^+$ for some
infinite cardinal $\kappa$. Let $(A_\alpha)_{\alpha < \kappa^+}$ be an increasing
continuous sequence of subsets of $V$ such that $|A_\alpha| < \kappa^+$ for all
$\alpha < \kappa^+$ and $V = \bigcup\{A_\alpha\colon \alpha < \kappa^+\}$.
For $\alpha < \kappa^+$ we define
\begin{eqnarray*}
V_\alpha & := & (V\setminus A_\alpha) \cup f^{-1}(\kappa^+)\\
E_\alpha & := & \big\{\{x,y\}\in E\colon x\in V_\alpha, \,y\in V\setminus A_\alpha \big\}\\
G_\alpha & := & (V_\alpha, E_\alpha)\\
f_\alpha & := & f \upharpoonright V_\alpha
\end{eqnarray*}
For any property $P$, $(A_\alpha)_{\alpha < \kappa^+}$ is said to be a 
\begr{$P$-destruction} of $(G,f)$ if $$S=\{\alpha < \kappa^+\colon (G_\alpha, f_\alpha)
\text{ does not fulfill }P\}$$ is stationary in $\kappa^+$. $(G,f)$ is called 
\begr{$P$-destructed} if there is a $P$-destruction of $(G,f)$.

\paragraph{Lemma~4} \textit{(Transfer Lemma)} 
Let $P(G,f)$ be a necessary condition for the existence of a perfect
$f$-factor of a graph $G$. If $(G,f)\in\cal C$, $|V(G)|=\kappa^+$ for an infinite
cardinal $\kappa$, and if $G$ possesses a perfect $f$-factor, then $(G,f)$ is not $P$-destructed.

\paragraph{Proof} Let $F$ be a perfect $f$-factor of $G$ and assume
that there is a $P$-destruction $(A_\alpha)_{\alpha<\kappa^+}$ of $(G,f)$. Define $V_\alpha,
E_\alpha, G_\alpha, f_\alpha, S$ as above and let $\alpha \in S$. $(G_\alpha, f_\alpha)$
does not fulfill $P$, and by the hypothesis of the Lemma, $G_\alpha$ has not a perfect
$f_\alpha$-factor. In particular $F_\alpha := F\cap E_\alpha$ is not a perfect $f_\alpha$-factor
of $G_\alpha$. Therefore there is a vertex $x_\alpha\in V_\alpha$ such that
$d_{F_\alpha}(x_\alpha)<f_\alpha(x_\alpha)=f(x_\alpha)$. Since $F$ is a perfect $f$-factor, there
exists, for some vertex $y_\alpha$, an edge $\{x_\alpha, y_\alpha\}\in F\setminus F_\alpha$. Using the fact
$|A_\alpha| < \kappa^+$ we know that $d_{F_\alpha}(x) = d_F(x) = \kappa^+=f(x)$ for any
$x\in f^{-1}(\kappa^+)$. So $x_\alpha\in V_\alpha\setminus f^{-1}(\kappa^+)$ and 
$y_\alpha\in A_\alpha\setminus f^{-1}(\kappa^+)$. .

If $\alpha \in S$ is a limit ordinal, let $\beta(\alpha) < \alpha$ be an
ordinal with $y_\alpha\in A_{\beta(\alpha)}$. By Fodor's Theorem
(cf.~\cite{1} or \cite{2}, Theorem~1.8.8), there is an ordinal 
$\gamma < \kappa^+$ such that
$$|\{\alpha\in S\colon \alpha\text{ limit ordinal, } \beta(\alpha) = \gamma\}| = \kappa^+.$$
Since $|A_\gamma| < \kappa^+$, there is a vertex $y^*\in A_\gamma$ such that
$$|\{\alpha\in S\colon \alpha\text{ limit ordinal, } y_\alpha = y^*\}| = \kappa^+.$$

If $x\in A_{\alpha_0}\setminus f^{-1}(\kappa^+)$ for some $\alpha_0<\kappa^+$,
then $x\not\in V_\alpha$ for all $\alpha > \alpha_0$
and thus
$$|\{\alpha \in S\colon x_\alpha = x\}| 
 < \kappa^+.$$
It follows that $f(y^*) = d_F(y^*) = \kappa^+$, so 
$y^*\in f^{-1}(\kappa^+)$. On the other hand
$y^*\in A_\alpha\setminus f^{-1}(\kappa^+)$
for every ordinal $\alpha$ with $y^*=y_\alpha$. This
contradiction proves the lemma.

\paragraph{Theorem~2}
Let $(G,f)\in\cal C$ and $|V(G)| = \aleph_1$. If $P$ is a hereditary
property such that $P(G,f)$ and if $(G,f)$ is not $P$-destructed then $G$
possesses a perfect $f$-factor.

\paragraph{Proof} Let $(A_\alpha)_{\alpha < \womega_1}$ be an increasing continuous
sequence of countable subsets of $V$ such that $V=\bigcup_{\alpha < \womega_1} A_\alpha$.
Define $V_\alpha, E_\alpha, G_\alpha, f_\alpha$ as above. Since 
$(A_\alpha)_{\alpha<\womega_1}$ is not a $P$-destruction, there is a closed unbounded 
set $K \subseteq \womega_1$ such that $(G_\alpha, f_\alpha)$ fulfills $P$ for every 
$\alpha\in K$. We can assume w.\,l.\,o.\,g.\ that $K=\womega_1$, because otherwise we could 
consider the sequence $(A_\alpha)_{\alpha\in K}$ instead of $(A_\alpha)_{\alpha< 
\womega_1}$. Since $(G,f)$ fulfills $P$ we can further assume that $A_0 = \emptyset$.

To obtain a perfect $f$-factor of $G$, we now construct an increasing continuous function
\mbox{$i\colon \womega_1 \rightarrow \womega_1$} and an increasing sequence 
$(F_\varepsilon)_{\varepsilon< \womega_1}$ of $f$-factors of $G$ with the following 
properties:
\begin{rom}
\item $\bigcup F_\varepsilon \subseteq A_{i(\varepsilon)}$
\item $\forall x\in A_{i(\varepsilon)}\big(f(x) \leq \aleph_0 \Rightarrow 
d_{F_\varepsilon}(x) = f(x)\big)$
\item $\forall x\in A_{i(\varepsilon)}\big(f(x) = \aleph_1 \Rightarrow
d_{F_{\varepsilon+1}\setminus F_\varepsilon} (x)=\aleph_0\big)$
\end{rom}
Then $F:= \bigcup_{\varepsilon < \womega_1} F_\varepsilon$ obviously is a perfect $f$-factor
of $G$.

The function $i$ and the sequence $(F_\varepsilon)_{\varepsilon < \womega_1}$ will be
defined by transfinite recursion. Let $i(0) := 0$ and $F_0 := \emptyset$. Now let
$\varepsilon >0 $ and let us assume that, for each $\delta < \varepsilon$, $i(\delta)$
and $F_\delta$ are already defined. If $\varepsilon$ is a limit ordinal, let
$i(\varepsilon) := \bigcup_{\delta < \varepsilon} i(\delta)$ and $F_\varepsilon 
:= \bigcup_{\delta<\varepsilon} F_\delta$.

Now let $\varepsilon = \delta+1$. By induction on $m$ we define an increasing sequence
$(H_m)_{m<\womega}$  of finite $f_{i(\delta)}$-factors of $G_{i(\delta)}$, an increasing
function $\varrho\colon\womega\longrightarrow\womega_1$, and, for any $n\geq m$, vertices 
$x_{m,n} \in V_{i(\delta)}$ such that for every $m$
\begin{enumerate}
\renewcommand{\theenumi}{(\alph{enumi})}
\renewcommand{\labelenumi}{\theenumi}
\item $\{x_{m,n}\colon n\geq m\} = A_{\varrho(m+1)}\setminus (A_{\varrho(m)}\setminus f^{-1}(\aleph_1))$
\item $\bigcup H_m \subseteq A_{\varrho(m)}$
\item $d_{H_{m+1}}(x_{k,m})=f_{i(\delta)}(x_{k,m})$ for all $k\leq m$ with $f_{i(\delta)}(x_{k,m})<\aleph_0$
\item $d_{H_{m+1}\setminus H_m}(x_{k,m})>0$ for all $k\leq m$ with $f_{i(\delta)}(x_{k,m})\geq\aleph_0$
\item $P(G_{i(\delta)} - H_m, f_{i(\delta)} - d_{H_m})$.
\end{enumerate}
Then let $F_\varepsilon := F_\delta\cup\bigcup\{H_m\colon m<\womega\}$ and
$i(\varepsilon):=\bigcup\{\varrho(m)\colon m<\womega\}$. By construction,
(i), (ii), (iii) are fulfilled.
\begin{entry}
\item[$\mathbf{m=0\colon}$]
Let $\varrho(0):=i(\delta)$, $H_0:=\emptyset$.

\item[$\mathbf{m=m+1\colon}$] Now suppose that for $m<\womega$ the ordinal
$\varrho(m)$, the finite $f_{i(\delta)}$-factor $H_m$ of $G_{i(\delta)}$, and, 
for all $k<m$ and $n\geq k$,
the vertices $x_{k,n}\in V_{i(\delta)}$ are already defined such that (a) - (e)
are fulfilled.

The set $W_m:=\{x_{k,n}\colon k\leq n<m\}$ is finite. Since $P$ is hereditary,
there exists a finite 
$f_{i(\delta)}$-factor $H_{m+1}\supseteq H_m$ of $G_{i(\delta)}$ such that 
$P(G_{i(\delta)} - H_{m+1}, f_{i(\delta)} - d_{H_{m+1}})$ and
$d_{H_{m+1}}(x) = f_{i(\delta)}(x)$ whenever $x\in W_m$ and
$f_{i(\delta)}(x) < \aleph_0$, or $d_{H_{m+1}\setminus H_m}(x) > 0$
whenever $x\in W_m$ and $f_{i(\delta)}(x) \geq \aleph_0$.

Let $\varrho(m+1) > \varrho(m)$ be the least ordinal such that
$\bigcup H_{m+1} \subseteq A_{\varrho(m+1)}$. For $n\geq m$ choose $x_{m,n}$ with 
$\{x_{m,m}, x_{m,m+1}, x_{m,m+2}, \ldots\} = A_{\varrho(m+1)}\setminus 
(A_{\varrho(m)}\setminus f^{-1}(\aleph_1))$.
\end{entry}

\paragraph{Corollary~2} Let $(G,f)\in\cal C$ and $|V(G)| = \aleph_1$.
\begin{rom}
\item
$G$ possesses a perfect $f$-factor if and only if $(G, f)$ is not $P_2$-destructed.
\item
If $ran(f)\subseteq\womega$
then $G$ possesses a perfect
$f$-factor if and only if $(G,f)$ is not $P_4$-destructed.
\end{rom}

\bigskip

To handle the cases of higher cardinality, we introduce the notion of a $\kappa$-perfect 
$f$-factor.

\paragraph{Definition} 
Let $(G,f)\in\cal C$ and let $\kappa$ be an infinite cardinal.
An $f$-factor $F$ of $G$ is said to be \begr{$\kappa$-perfect} if $d_F(x) = f(x)$
for all vertices $x$ with $f(x) \leq \kappa$ and
$d_F(x)>0$ for all vertices $x$ with $f(x)>\kappa$.

\paragraph{Theorem~3}
Let $\kappa$ be an infinite cardinal, $(G,f)\in\cal C$, and $|V(G)| = \kappa^+$.
$G$ possesses a perfect $f$-factor if and only if there is an increasing
continuous sequence $(A_\alpha)_{\alpha<\kappa^+}$ of subsets of $V(G)$ such that
\begin{rom}
\item $A_0 = \emptyset$, $V(G) = \bigcup\{A_\alpha\colon \alpha < \kappa^+\}$,
\item $|A_{\alpha+1} \setminus A_\alpha| = \kappa$ for all $\alpha < \kappa^+$,
\item
for all $\alpha < \kappa^+$ there exists an $\kappa$-perfect $g_\alpha$-factor
of\\
 $(B_\alpha,\big\{\{x,y\}\in E\colon x\in B_{\alpha}, y\in A_{\alpha+1} \setminus 
A_\alpha\big\})$, where\\
$B_\alpha=(A_{\alpha+1} \setminus (A_\alpha \setminus f^{-1}(\kappa^+)))$ and
$g_\alpha := f\upharpoonright B_\alpha$.
\end{rom}

\paragraph{Proof}
Let $(A_\alpha)_{\alpha<\kappa^+}$ be an increasing continuous
sequence of subsets of $V$ and, for $\alpha<\kappa^+$, let $F_\alpha$ be a
$\kappa$-perfect $g_\alpha$-factor with the properties (i), (ii), (iii).
Then $F_{\alpha_1}\cap F_{\alpha_2}=\emptyset$ if $\alpha_1\neq\alpha_2$.
Let $F:=\bigcup\{F_\alpha\colon\alpha<\kappa^+\}$.
We will show that $F$ is a perfect $f$-factor of $G$. Let $x\in V$ and let $\alpha$
be the smallest ordinal such that $x\in A_{\alpha+1}$. If $f(x)\leq\kappa$ then
$d_F(x)=d_{F_\alpha}(x)=f(x)$. If on the other hand $f(x)>\kappa$, we have
$d_{F_\beta}(x)>0$ for all $\beta\geq \alpha$ since $F_\beta$ is $\kappa$-perfect.
Thus $d_F(x)=\kappa^+$.

To prove the converse, let $F$ be a perfect $f$-factor of $G$ and $A_0:=\emptyset$.
Let $(P_{\delta}\colon\delta<\kappa^+)$ be a partition of $V$ such that
$|P_{\delta}| = \kappa$ for all $\delta < \kappa^+$.
Now assume that $A_\delta \subseteq V$ is defined for all $\delta < \alpha$.
If $\alpha$ is a limit ordinal, then let $A_\alpha=\bigcup\{A_\delta\colon\delta<\alpha\}$.
If $\alpha=\delta+1$, we define by induction an increasing sequence $(C_n)_{n<\womega}$
of subsets of $V$. Let $C_0\subseteq V$ such that
$A_\delta\cup P_\delta\subseteq C_0$ and $|C_0\setminus A_\delta| = \kappa$. If $C_n$ is
defined let $C_{n+1}$ be a ''$\kappa$-neighborhood'' of $C_n$: If $x\in C_n$ and
$f(x)\leq\kappa$ let $N(x)=\{y\in V\colon \{y,x\}\in F\}$, and if $f(x)=\kappa^+$ choose
$y_x\in V\setminus C_n$ with $\{y_x,x\}\in F$ and let $N(x)=\{y_x\}$. Then let
$C_{n+1}=C_n\cup\bigcup\{N(x)\colon x\in C_n\}$ and
$A_{\alpha} := \bigcup\{C_n\colon n<\womega\}$.  By construction, 
$(A_\alpha)_{\alpha<\kappa^+}$ is an increasing continuous sequence of subsets of $V$ 
with the properties (i), (ii), (iii).

\paragraph{Remark}
If $\kappa^+=\aleph_2$, $g_\alpha := f \upharpoonright V_\alpha \cap A_{\alpha+1}$ and $X_\alpha :=
A_{\alpha+1} \cap f^{-1}(\aleph_2)$, then there is an $\aleph_1$-perfect $g_\alpha$-%
factor of $(V_\alpha\cap A_{\alpha+1}, \{\{x,y\}\in E\colon x\in V_\alpha \cap 
A_{\alpha+1}, y\in A_{\alpha+1}\setminus A_\alpha\})$ if and only if there exists a 
function $h_\alpha\colon A_{\alpha+1} \cap f^{-1}(\aleph_2) \rightarrow \womega\cup 
\{\aleph_0, \aleph_1\}$ such that there is a perfect $(g_\alpha\setminus (g_\alpha\!\! 
\upharpoonright
\!\! X_\alpha)) \cup h_\alpha$-factor of $(V_\alpha \cap A_{\alpha+1}, \{ \{x,y\}\in E \colon
x\in V_\alpha \cap A_{\alpha+1}, y\in A_{\alpha+1} \setminus A_\alpha\})$.

\paragraph{Corollary~3}
Let $(G,f)\in\cal C$ and $|V(G)| = \aleph_2$. $G$ possesses a perfect $f$-factor
if and only if there is an increasing continuous sequence $(A_\alpha)_{\alpha < 
\womega_2}$ of subsets of $V(G)$, such that 
\begin{rom}
\item $A_0 = \emptyset$, $V(G) = \bigcup_{\alpha < \womega_2} A_\alpha$.
\item $|A_{\alpha+1}\setminus A_\alpha| = \aleph_1$ for all $\alpha < \womega_2$.
\item For each $\alpha < \womega_2$ there is a function $h_\alpha\colon A_{\alpha+1} \cap
f^{-1}(\aleph_2) \rightarrow \womega\cup \{\aleph_0, \aleph_1\}$ such that the graph
$$(A_{\alpha+1}\setminus (A_\alpha \setminus f^{-1}(\aleph_2)), \{\{x,y\}\in E\colon
x\in A_{\alpha+1} \setminus (A_\alpha \setminus f^{-1}(\aleph_2)), y\in A_{\alpha+1}
\setminus A_\alpha\})$$
together with 
$$\big(f \upharpoonright
A_{\alpha+1} \setminus (A_\alpha \setminus f^{-1}(\aleph_2)\,) \, \setminus \,
f \upharpoonright (A_{\alpha+1} \cap f^{-1}(\aleph_2))\,\big) \cup h_\alpha$$
is not $P_2$-destructed.
\end{rom}

\end{document}